# CANONICAL MOMENTA OF NONLINEAR COMBAT

LTC Michael Bowman, USA
*U.S. Army Materiel Command, ATTN: AMCRDA-TF, 5001 Eisenhower Avenue Alexandria, VA 22333-0001*
e-mail: mbowman@hqamc.army.mil, mbowman3@osf1.gmu.edu
and
Lester Ingber
*Lester Ingber Research, P.O. Box 857, McLean, VA 22101*
e-mail: ingber@ingber.com, ingber@alumni.caltech.edu

**ABSTRACT**

The context of nonlinear combat calls for more sophisticated measures of effectiveness. We present a set of tools that can be used as such supplemental indicators, based on stochastic nonlinear multivariate modeling used to benchmark Janus simulation to exercise data from the U.S. Army National Training Center (NTC). As a prototype study, a strong global optimization tool, adaptive simulated annealing (ASA), is used to explicitly fit Janus data, deriving coefficients of relative measures of effectiveness, and developing a sound intuitive graphical decision aid, canonical momentum indicators (CMI), faithful to the sophisticated algebraic model. We argue that these tools will become increasingly important to aid simulation studies of the importance of maneuver in combat in the 21st century.

## 1. INTRODUCTION

Add to the fog of war, the haze of information overload. As we attempt to use technology to remove the fog of war from military operations we face new challenges in sifting through potential mountains of available information and new battlefield capabilities to make more complex decisions in a more rapid manner. The correct tools can help, using no tools or the wrong tools can be disastrous.

### 1.1. Characteristics of Modern Combat Operations — The Non-Linear Battlefield

The U.S. Army Training and Doctrine Command (TRADOC) is working to develop operational concepts for land combat in the 21st century (U.S. Army Training and Doctrine Command, 1996). TRADOC characterizes the expected operational environment of the future as: multi-dimensional, precise, and non-linear, with distributed operations and simultaneity.

> Multi-dimensional because beyond the traditional width, depth and height, the factors of time, the electromagnetic spectrum and the human dimension (soldiers, leaders and civilian populations) must be accounted for.
>
> Precise because friendly and enemy targets will be attacked with precision by smart and even "brilliant" weapons throughout the battlespace.
>
> Non-linear because in the entire spectrum of operational environments, peacekeeping to high intensity large scale warfare, military operations will involve the accomplishment of tasks across the entire battlespace rather than straightforward massing of combat power along the traditional forward line of troops (FLOT).
>
> With distributed operations because effective operations will be conducted throughout the battlespace, where and when required to achieve decisive effects vice concentrated at a single, possibly decisive point.
>
> With simultaneity because as the other characteristics are put into motion in numerous, simultaneous operations the enemy will be presented with multiple crisis with few options for an effective response.

These characteristics are not necessarily new or poorly understood. Combat has almost always been known to be a complex, confusing and horrific endeavor. Successful commanders in documented battles and military operations from the days of Sun Tzu (Tzu, 1963) to Operations Just Cause and Desert Storm have attempted to use these factors to their advantage. The attempts to capture and study the theories war are well documented in many works (Leonhard, 1991).

What is revolutionary is that the U.S. Army is seeking to systematically put the advantages of these characteristics to work in military operations while mitigating the risks posed by them. The Army hopes to accomplish this through the development and use of modernized doctrine, tactics-techniques and procedures (TTP), as well the introduction of advanced automation and communications equipment commonly call "digitization."

TRADOC defines digitization of the battlespace as "the application of technology to acquire, exchange, and employ timely information horizontally and vertically integrated to create a common picture of the battlefield from soldier to commander." Thus digitization attempts to lift some of the fog of war through the concerted use of information. Digitization is a key piece of the Army's efforts to move the force more fully into what has been called the "next wave" of warfare (Toffler and Toffler, 1993) or "information warfare."

### 1.2. Analysis of the Nonlinear Battlefield

Too often the management of complex systems is ill-served by not utilizing the best tools available. For example, requirements set by decision-makers often are not formulated in the same language as constructs formulated by powerful mathematical formalisms, and so the products of analyses are not properly or maximally utilized, even if and when they come close to faithfully representing the powerful intuitions they are supposed to model. In turn, even powerful mathematical constructs are ill-served, especially when dealing with multivariate nonlinear complex systems, when these formalisms are butchered into quasi-linear approximations to satisfy constraints of numerical algorithms familiar to particular analysts, but which tend to destroy the power of the intuitive constructs developed by decision-makers. These problems are present in many disciplines.

For at least a large class of systems, including some classes of large-scale combat, these problems can be bypassed by using a blend of an intuitive and powerful mathematical-physics formalism to generate "canonical momenta" indicators (CMI), which are used by AI-type rule-based models of management of complex systems. Typically, both the formalism generating the CMI and the rule-based models have quite nonlinear constructs, and they must be "trained" or fit to data subsequent to testing on "out-of-sample" data, before they can be used effectively for "real-time" production runs. To handle these fits of nonlinear models of real-world data, a generic powerful optimization code, Adaptive Simulated Annealing (ASA), has been developed (Ingber, 1993a).

The algebraic and numerical methodology closely follows modeling recently published by one of us (LI) (Ingber, 1997a) in analyses of electroencephalography (EEG) (Ingber, 1997b) and finance (Ingber, 1996b; Ingber, 1996c).

## 2. BACKGROUND

### 2.1. The U.S. Army National Training Center (NTC)

The NTC is a large maneuver range dedicated to the simulation of desert combat, training battalion and brigade size mechanized units from U.S. Army heavy divisions and separate brigades. The NTC is unique in that it is highly instrumented with the Multiple Integrated Laser Engagement System (MILES) and range instrumentation which follows the location and activity of most vehicles and some dismounted infantry. The NTC also has a dedicated Opposing Force (OPFOR) which acts as the enemy during force-on-force exercises with visiting units.

Transfers of data between different databases and computer operating systems were automated by one of us (MB) (Bowman, 1989), He has coordinated and integrated data from NTC, Training and Doctrine Command (TRADOC) Analysis Command (TRAC) at White Sands Missile Range, New Mexico (TRAC-WSMR) and at Monterey, California (TRAC-MTRY) for Janus(T) wargaming at TRAC-MTRY, and for use at at Lawrence Livermore National Laboratory (LLNL) Division B, and for Janus(T) and NTC modeling.

### 2.2. Janus

Janus is an interactive, two-sided, closed, stochastic, ground combat simulation. Players direct their elements, executing tactical plans and reacting to enemy actions. The disposition of opposing forces is not completely known to players. Janus models individual systems moving, searching, detecting and engaging other ground or air systems over a three-dimensional terrain representation, using Army-developed algorithms and data to model combat processes.

### 2.3. Statistical Mechanics of Combat

A series of papers has developed a statistical mechanics of large-scale combat (Ingber, 1985; Ingber, 1986; Ingber, 1993c; Ingber, Fujio, and Wehner, 1991; Ingber and Sworder, 1991), where details and the rationale of this presentation can be found. The statistical mechanics of combat (SMC) modeling approach used here was developed by LI when he was principal investigator of an Army contract to benchmark Janus simulation to NTC exercise data.

Consider a scenario taken from our NTC study: two Red systems, Red T-72 tanks ($RT$) and Red armored personnel carriers ($RBMP$), and three Blue systems, Blue M1A1 and M60 tanks ($BT$), Blue armored personnel carriers ($BAPC$), and Blue tube-launched optically-tracked wire-guided missiles ($BTOW$), where $RT$ specifies the number of Red tanks at a given time $t$, etc. Consider the kills suffered by $BT$, $\Delta BT$, e.g., within a time epoch $\Delta t \approx 5$ min

$$\Delta BT/\Delta t \equiv \dot{BT} = x^{BT}_{RT} RT + y^{BT}_{RT} RT\ BT + x^{BT}_{RBMP} RBMP + y^{BT}_{RBMP} RBMP\ BT \tag{1}$$

Here, the $x$ terms represent attrition owing to point fire; the $y$ terms represent attrition owing to area fire. Note that the algebraic forms chosen are consistent with current perceptions of aggregated large scale combat. The version of Janus(T) used to generate this data does not permit direct-fire fratricide; such terms are set to zero. In most NTC scenarios fratricide typically is negligible.

Now consider sources of noise, e.g., that at least arise from PD, PA, PH, PK, etc. Furthermore, such noise likely has its own functional dependencies, e.g., possibly being proportional to the numbers of units involved in the combat. For simplicity here, still generating much nonlinearity, only diagonal noise terms are considered. Coupling among the variables takes place in the drift terms (deterministic limit); for simplicity only linear terms in the drifts are taken for this prototype study.

$$\frac{\Delta BT}{\Delta t} \equiv \dot{BT} = x^{BT}_{RT} RT + y^{BT}_{RBMP} RT\ BT + x^{BT}_{RBMP} RBMP + y^{BT}_{RBMP} RBMP\ BT + z^{BT}_{BT} BT \eta^{BT}_{BT} \tag{2}$$

where the $\eta$ represent sources of (white) noise (in the Itô prepoint discretization discussed below). The noise terms are taken to be log normal (multiplicative) noise for the diagonal terms and additive noise for the off-diagonal terms. This induces a high degree of nonlinearity, which can be seen by transforming each variable $M^G$ to $X^G$

$$X^G = \ln M^G\ ,\ \ \dot{X}^G = \dot{M}^G / M^G\ ,$$

$$M^G = \{RT,\ RBMP,\ BT,\ BAPC,\ BTOW\} \tag{3}$$

yielding $X^G$ equations with constant coefficients of the noise, at the expense of introducing exponential terms in the drifts.

The methodology presented here can accommodate any other nonlinear functional forms, and any other variables that can be reasonably represented by such rate equations, e.g., expenditures of ammunition or bytes of communication (Ingber, 1989a). Variables that cannot be so represented, e.g., terrain, $C^3$, weather, etc., must be considered as "super-variables" that specify the overall context for the above set of rate equations.

### 2.4. Janus Data

For this study, data collected during our NTC-Janus(T) project circa 1988 was used to fit the coefficients of the above 5 coupled equations. Time epochs were 5 mins each, and we used data from 6 battle simulations between 30 mins and 75 mins into the battles, for a total of 60 states of data, each state giving the present values of each of the 2 Red and 3 Blue units.

It should be noted that the numbers of units in this particular set of data are barely large enough to be considered large-scale so that the statistical methodology being presented is applicable. At the least, this paper presents a full study to demonstrate the SMC approach for future sets of large-scale data.

### 2.5. Algebraic Development

The five coupled stochastic differential equations, for variables $M^G = \{RT, RBMP, BT, BAPC, BTOW\}$, can be represented equivalently by a short-time conditional probability distribution $P$ in terms of a Lagrangian $L$:

$$P(R\cdot, B\cdot; t + \Delta t | R\cdot, B\cdot; t) = \frac{1}{(2\pi\Delta t)^{5/2}\sigma^{1/2}} \exp(-L\Delta t) \quad (4)$$

where $\sigma$ is the determinant of the inverse of the covariance matrix, the metric matrix of this space, "$R\cdot$" represents $\{RT, RBMP\}$, and "$B\cdot$" represents $\{BT, BAPC, BTOW\}$. (Here, the prepoint discretization is used, which hides the Riemannian corrections explicit in the midpoint discretized Feynman Lagrangian; only the latter representation possesses a variational principle useful for arbitrary noise.)

This defines a scalar "dynamic cost function," $C(x, y, z)$,

$$C(x, y, z) = L\Delta t + \frac{5}{2}\ln(2\pi\Delta t) + \frac{1}{2}\ln\sigma \quad (5)$$

which can be used with the adaptive simulated annealing (ASA) algorithm (Ingber, 1989b; Ingber, 1993a) further discussed below, to find the (statistically) best fit of $\{x, y, z\}$ to the data.

The form for the Lagrangian $L$ and the determinant of the metric $\sigma$ to be used for the cost function $C$ is

$$L = \sum_G \sum_{G'} \frac{(\dot{M}^G - g^G)(\dot{M}^{G'} - g^{G'})}{2g^{GG'}}$$

$$\sigma = \det(g_{GG'})$$

$$(g_{GG'}) = (g^{GG'})^{-1}$$

$$g^{GG'} = \sum_i \hat{g}_i^G \hat{g}_i^{G'} \quad (6)$$

It must be emphasized that the output need not be confined to complex algebraic forms or tables of numbers. Because $\underline{L}_F$ possesses a variational principle, sets of contour graphs, at different long-time epochs of the path-integral of $P$, integrated over all its variables at all intermediate times, give a visually intuitive and accurate decision aid to view the dynamic evolution of the scenario. Also, this Lagrangian approach permits a quantitative assessment of concepts usually only loosely defined, which are used to advantage here.

$$\text{Momentum} \equiv \Pi^G = \frac{\partial \underline{L}_F}{\partial(\partial M^G/\partial t)} \ ,$$

$$\text{Mass} \equiv g_{GG'} = \frac{\partial \underline{L}_F}{\partial(\partial M^G/\partial t)\partial(\partial M^{G'}/\partial t)} \ ,$$

$$\text{Force} \equiv \frac{\partial \underline{L}_F}{\partial M^G} \ ,$$

$$F - ma \equiv \delta \underline{L}_F = 0 = \frac{\partial \underline{L}_F}{\partial M^G} - \frac{\partial}{\partial t}\frac{\partial \underline{L}_F}{\partial(\partial M^G/\partial t)} \ . \quad (7)$$

These momenta are the canonical momenta indicators (CMI).

### 2.6. Numerical Methods

A systematic numerical procedure has been developed for fitting parameters in such stochastic nonlinear systems to data using methods of adaptive simulated annealing (ASA) as maximum likelihood technique on the Lagrangian (Ingber, 1989b; Ingber, 1993a; Ingber, 1993b; Ingber, 1996a), and then integrating the path integral using a non-Monte Carlo technique especially suited for nonlinear systems (Wehner and Wolfer, 1983). This numerical methodology has been applied with success to several systems (Ingber, 1990; Ingber, 1991; Ingber, 1995; Ingber, 1996b; Ingber, Fujio, and Wehner, 1991; Ingber and Nunez, 1990; Ingber and Nunez, 1995). ASA has been applied to many problems by many people in many disciplines (Ingber, 1993b; Ingber, 1996a).

The feedback of many users regularly scrutinizing the source code ensures its soundness as it becomes more flexible and powerful. The ASA code can be obtained at no charge, via WWW from http://www.ingber.com/, or via FTP from ftp.ingber.com. The file http://www.ingber.com/MISC.DIR/asa_examples has several templates of "toy" test problems, especially illustrating how tuning can increase the efficiency of ASA by orders of magnitude. The file http://www.ingber.com/asa_papers has references to the the use of ASA by some other researchers, e.g., in studies ranging from: comparisons among SA algorithms and between ASA and genetic algorithms, tabu and hillclimbing (Ingber and Rosen, 1992; Mayer *et al*, 1996; Rosen, 1992), to molecular models (Su *et al*, 1996), to imaging (Wu and Levine, 1993), to neural networks (Cohen, 1994), to econometrics (Sakata, 1995), to geophysical inversion (Sen and Stoffa, 1995), to wide-spread use in financial institutions (Wofsey, 1993), etc.

## 3. PRESENT RESULTS

### 3.1. Janus

Table 1 gives the results of ASA fits of the above 5 coupled equations to Janus-generated data. Note that the noise coefficient is roughly the same for all units, being largest for *BTOW*. Note the relative importance of coefficients in "predicting" the immediate next epoch, with *BTOW* larger than *BAPC* larger than *BT* in depleting Red forces (but being multiplied by the total number of units at any time). The coefficients of "prediction" of attrition by Red forces has *RT* larger than *RBMP* against *BTOW*, and *RT* less than *RBMP* against *BT* and *BAPC* (but being multiplied by the total number of units at any time).

| | *RT* | *RBMP* | *BT* | *BPAC* | *BTOW* | $\eta$ [.] |
|---|---|---|---|---|---|---|
| $\dot{RT}$ | - | - | -8.6E-5 | -5.9E-3 | -3.6E-2 | 3.7E-3 |
| $\dot{RBMP}$ | - | - | -2.7E-3 | -2.2E-2 | -3.1E-2 | 4.3E-3 |
| $\dot{BT}$ | -6.7E-4 | -4.7E-3 | - | - | - | 7.9E-3 |
| $\dot{BAPC}$ | -1.0E-4 | -4.0E-3 | - | - | - | 6.7E-3 |
| $\dot{BTOW}$ | -2.1E-3 | -1.2E-6 | - | - | - | 1.3E-2 |

TABLE 1. Entities in the table are the ASA-fitted coefficients of the coupled set of 5 equations representing the dynamics of Red and Blue interactions. Note that the last column coefficients are multiplied by the corresponding variable in the first column. A dash represents no coefficient present in the equations.

The upper graph in Figure 1 gives the attrition data. The attrition data is given as the overage over 6 runs for each time point. The lower figure in Figure 1 gives the derived CMI. After the ASA fits, the CMI are calculated for each point in time in each of the 6 runs. The figure gives the average over the 6 runs for each time point. Note that the attrition rate of all units is fairly constant, and so there are no surprises expected in this kind of analysis. The marked changes of the systems at the end of the epoch signals the essential ending of the combat.

Using the particular model considered here, the CMI are seen to be complementary to the attrition rates, being somewhat more sensitive to changes in the battle than the raw data. The coefficients fit to the combat data are modifiable to fit the current "reality" of system capabilities.

### 3.2. Statistical Mechanics of Neocortical Interactions (SMNI)

In the context of this present study, the CMI are more sensitive measures than the energy density, effectively the square of the CMI, or the information which also effectively is in terms of the square of the CMI (essentially integrals over quantities proportional to the energy times a factor of an exponential including the energy as an argument). This is even more important when replenishment of forces is permitted, often leading to oscillatory variables. Neither the energy or the information give details of the components as do the CMI. The information and energy densities are calculated and printed out after fits to data, along with the CMI.

The utility of the CMI in such a system can be seen in Figure 2, from a recent study fitting SMNI to EEG data (Ingber, 1997b).

## 4. CONTEXT OF PRESENT STUDY

### 4.1. Janus Update

Not only are we moving to to a new era in tactics and doctrine with the theory of a nonlinear battlefield, and the "next wave" of warfare (information warfare), we've also seen a complete turn around in the capabilities of "Red" (old Soviet Union and client states) versus "Blue" (U.S./Nato) forces. When we did the studies of NTC and JANUS data in the late 1980's the Blue side was at a distinct technological disadvantage and the NTC scenarios were played out that way - the MILES sensors on T72s were positioned so that the T72s could not be killed by frontal hits by any U.S. weapons, while M60s could be killed by any hits from the T72 and just about anything on the battlefield could kill a U.S. APC or TOW vehicle.

In the 1990s the U.S. and NATO have advanced to a new generation of combat systems (M1A2 tank and Bradley Fighting Vehicle) while potential adversaries equipped with "Red" equipment (T72 and BMP) have not. This was dramatically apparent in the Gulf war in which M1s and Bradleys destroyed huge quantities of Iraqi equipment with almost no losses on the U.S. side. In fact the only M1 tanks destroyed in the gulf were hit by mistakenly by other M1 tanks. The U.S. Bradley fighting vehicle not only became a tank killer with its TOW missiles, it also killed everything short of tanks with its 30mm cannon. Thus the Bradley is a critical "killer" versus the "battlefield taxi" status the APC used by the U.S. in the 1980s.

The present study should be viewed as a prototype to similarly process new data as it becomes available.

### 4.2. Attrition Vs Maneuver Warfare

The "non-linear" battle field and the Army's modern "maneuver warfare" doctrines call for the a switch in emphasis from the fire and maneuver described in Air-Land Battle Doctrine, which carries a connotation of "attrition warfare" to an emphasis on the use of more pure maneuver to when ever possible by-pass and make irrelevant enemy strengths.

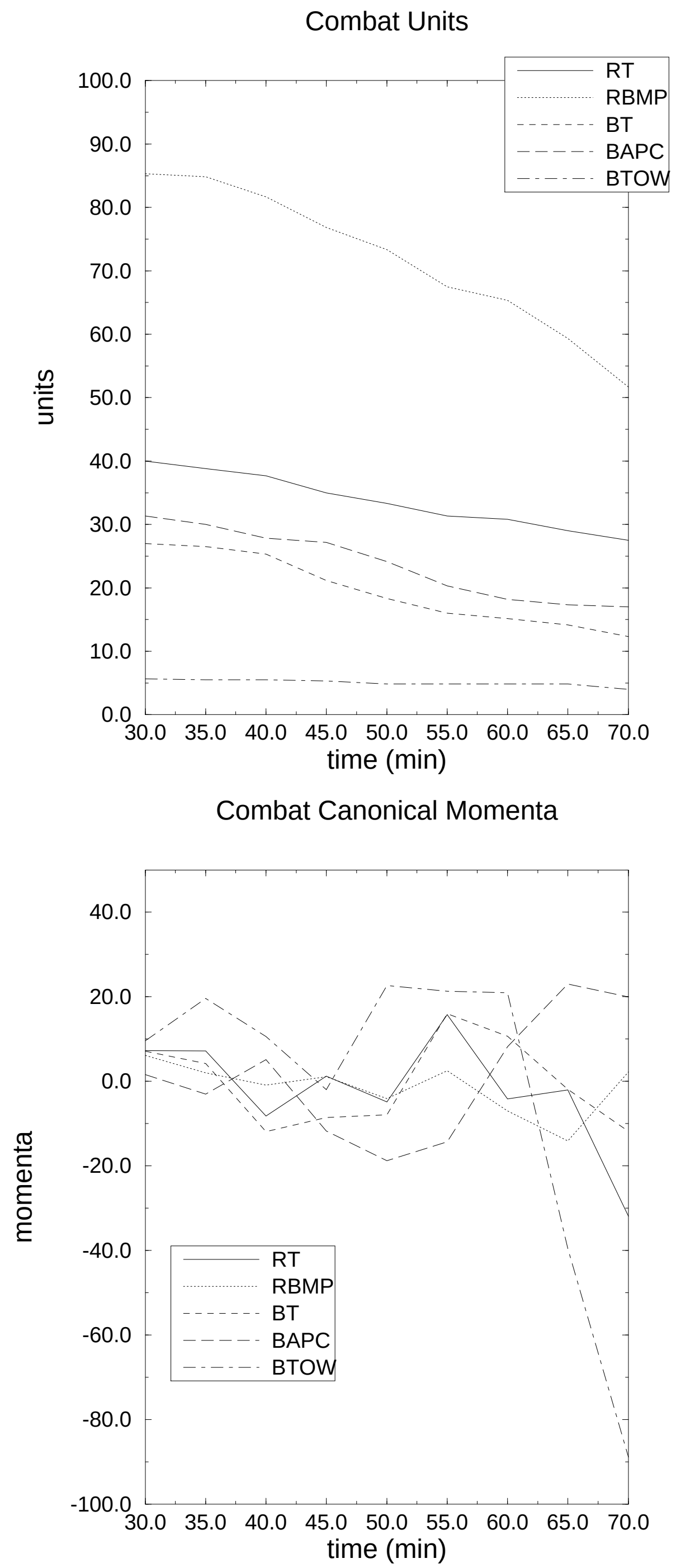

FIG. 1. The results of Janus(T) attrition of Red and Blue units are given in the upper figure. The canonical momenta indicators (CMI) for each system are given in the lower figure.

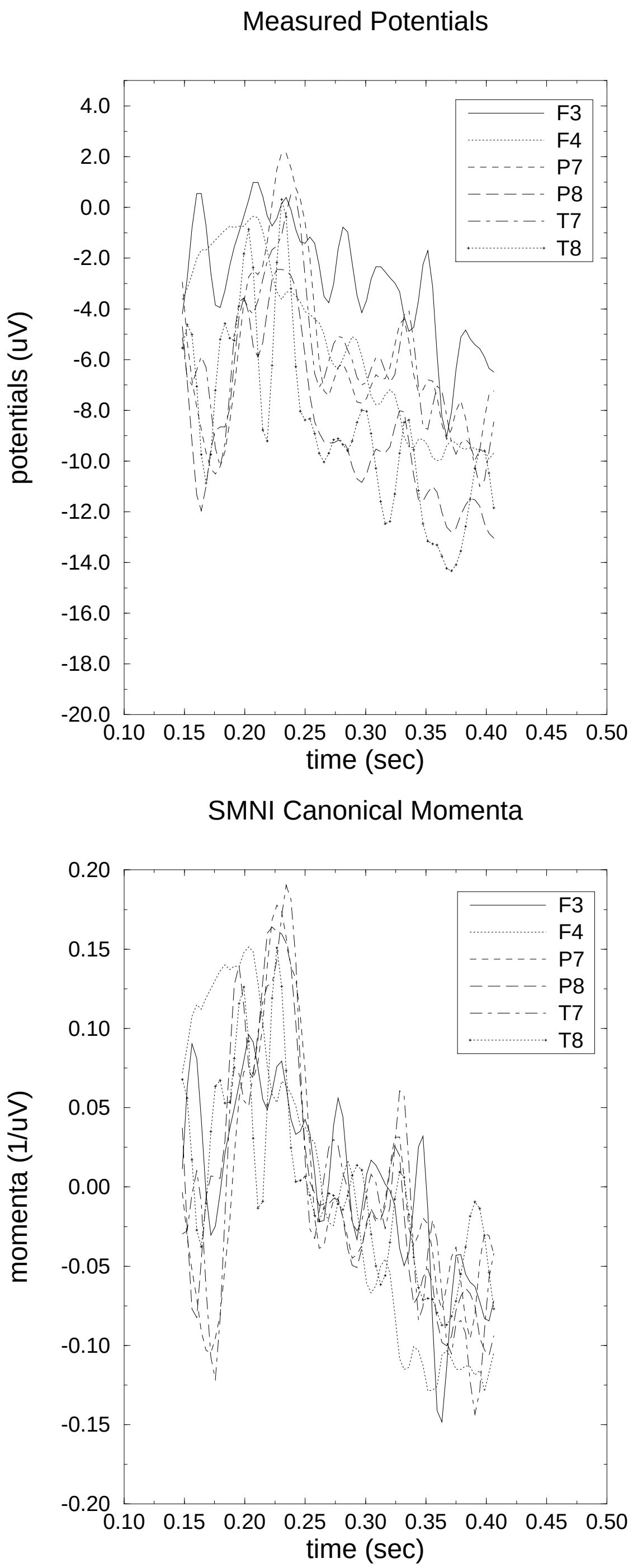

FIG. 2. For the matching paradigm for control subject co2c0000337, plots are given of activities under 6 electrodes of the electric potential in the top figure, and the CMI in the bottom figure.

The issue is that the strongest proponents of maneuver warfare may consider “force ratios”, “kill ratios”, “attrition rates,” etc., the tools of poor commanders who should be concentrating on finding and tipping the enemy center of gravity by maneuver rather than calculating ratios for head-on attacks.

While maneuver is the technique of choice, and driving the enemy from the battlefield without firing a shot is the goal, in all likelihood even the most masterfully maneuvered force will still fight engagements and battles during a campaign, and knowledge of ratios is still a valid command tool. Given the complex operation environment envisioned on a non-linear battlefield, some attrition combat is likely to be taking place at any given time that maneuver is being exercised on another portion of the battlefield.

In the context of the present study, the concept of attrition is still valid, even if one-sided. For example, given three possible strategies of maneuvers for a given upcoming battle, simulation studies can help address just which Blue units might be most effective in taking out various Red units. That is, in the context of the SMC papers, there really would not be any $d\mathrm{Blue}/dt$ equations (or they would be relatively insensitive in those context where there is little Blue attrition), just $d\mathrm{Red}/dt$ equations describing the attrition of Red forces due to various Blue forces, in the context of a given set of maneuvers. It can be argued that this is a necessary component of any planning, especially if and when reasonable measures of effectiveness are required for the nonlinear battlefield.

The approach presented here and in other SMC papers is more useful in the nonlinear battlefield than merely finding which Blue units take out which Red units, or vice versa. That is, in a nonlinear context, there is often an effective synergy among units of a force, such that a particular unit’s actual strength may not be measured in a very useful way just by correlating which opposite units it can attrite. The real measure of effectiveness is what the combined force can attrite on the opposing forces.

For example, a unit may be introduced as a measure of the communications network of a force. It is clear that statistical analyses of killer scoreboards will not suffice to measure the effectiveness of these units; the coefficients of the SMC equations can perform this function. In more technical terms, one must perform global optimization of the full multivariate stochastic system in order to reasonably measure the influence of any particular constituent.

## 5. CONCLUSIONS AND DIRECTIONS

We expect the CMI and the fitted coefficients to be more valuable predictors of events in combat, as the battlefield becomes more nonlinear. We have described a reasonable approach to quantitatively measuring this nonlinearity, and a reasonable approach to faithfully presenting this information to commanders in the field so that they may make timely decisions.

As was performed in the finance studies (Ingber, 1996b; Ingber, 1996c), a future project will similarly use recursive ASA optimization, with an inner-shell fitting CMI combat data, embedded in an outer-shell of parameterized customized commander’s AI-type rules acting on the CMI, to create supplemental decision aids.

Given the high cost of major field exercises in an environment of shrinking budgets, our forces will rely more and more heavily on modeling and simulation to develop, test, and practice tactics and doctrine at all levels. Modeling and simulation remain highly useful devices for making tactical mistakes and learning lessons at little or no cost. The use of CMI and ASA to evaluate and improve these models and simulation remains a worthy goal. This paper can be viewed as a call for data to perform future studies using this methodology.

The NTC remains possibly the best source of realistic, simulated combat data — as near to reality as we can get without bloodshed. The NTC instrumentation system has gone through one major upgrade since 1990 and another major upgrade is scheduled to be completed by the end of February 1997. This latest should allow for a greater variety of systems to be tracked, will more closely match firing systems and targets in the database, and will track up to 2000 systems. This upgrade is designed to support very detailed data collection on the Army's digitized Experimental Force (EXFOR) when it has its NTC rotation in Feb/March 1997. This NTC rotation is the culmination of the Army's Task Force XXI (TFXXI) Advanced Warfighting Experiment (AWE) which was designed and is being executed to test the effect of digitizing (providing shared battlefield awareness through computers and communications equipment) a Brigade sized task force. This upgrade to the NTC's digitization promises to make much more complete and accurate data available from the exercises done there.

**REFERENCES**

Bowman, M.. 1989. "Integration of the NTC Tactical Database and JANUS(T) Towards a Combat Decision Support System". M.S. Thesis. Naval Postgraduate School, Monterey, CA.

Cohen, B.. 1994. Training synaptic delays in a recurrent neural network, M.S. Thesis. Tel-Aviv University, Tel-Aviv, Israel.

Ingber, L.. 1985. "Statistical mechanics algorithm for response to targets (SMART)." In *Workshop on Uncertainty and Probability in Artificial Intelligence: UC Los Angeles, 14-16 August 1985*, American Association for Artificial Intelligence, Menlo Park, CA, pp. 258-264.

Ingber, L.. 1986. "Nonlinear nonequilibrium statistical mechanics approach to $C^3$ systems." In *9th MIT/ONR Workshop on $C^3$ Systems: Naval Postgraduate School, Monterey, CA, 2-5 June 1986*, MIT, Cambridge, MA, pp. 237-244.

Ingber, L.. 1989a. "Mathematical comparison of JANUS(T) simulation to National Training Center." In *The Science of Command and Control: Part II, Coping With Complexity*, S.E. Johnson and A.H. Levis, eds. AFCEA International, Washington, DC, pp. 165-176.

Ingber, L.. 1989b. "Very fast simulated re-annealing." *Mathl. Comput. Modelling 12* no. 8: 967-973.

Ingber, L.. 1990. "Statistical mechanical aids to calculating term structure models." *Phys. Rev. A 42* no. 12: 7057-7064.

Ingber, L.. 1991. "Statistical mechanics of neocortical interactions: A scaling paradigm applied to electroencephalography." *Phys. Rev. A 44* no. 6: 4017-4060.

Ingber, L.. 1993a. "Adaptive Simulated Annealing (ASA)". [http://www.ingber.com/ASA-shar, ASA-shar.Z, ASA.tar.Z, ASA.tar.gz, ASA.zip]. Lester Ingber Research, McLean, VA.

Ingber, L.. 1993b. "Simulated annealing: Practice versus theory." *Mathl. Comput. Modelling 18* no. 11: 29-57.

Ingber, L.. 1993c. "Statistical mechanics of combat and extensions." In *Toward a Science of Command, Control, and Communications*, C. Jones, eds. American Institute of Aeronautics and Astronautics, Washington, D.C., pp. 117-149.

Ingber, L.. 1995. "Path-integral evolution of multivariate systems with moderate noise." *Phys. Rev. E 51* no. 2: 1616-1619.

Ingber, L.. 1996a. "Adaptive simulated annealing (ASA): Lessons learned." *Control and Cybernetics 25* no. 1: 33-54.

Ingber, L.. 1996b. "Canonical momenta indicators of financial markets and neocortical EEG." In *International Conference on Neural Information Processing (ICONIP'96)*, Springer, New York, pp. 777-784.

Ingber, L.. 1996c. "Statistical mechanics of nonlinear nonequilibrium financial markets: Applications to optimized trading." *Mathl. Computer Modelling 23* no. 7: 101-121.

Ingber, L.. 1997a. "Data mining and knowledge discovery via statistical mechanics in nonlinear stochastic systems." : (submitted). .

Ingber, L.. 1997b. "Statistical mechanics of neocortical interactions: Applications of canonical momenta indicators to electroencephalography." *Phys. Rev. E* : (to be published).

Ingber, L., Fujio, H., and Wehner, M.F.. 1991. "Mathematical comparison of combat computer models to exercise data." *Mathl. Comput. Modelling 15* no. 1: 65-90.

Ingber, L. and Nunez, P.L.. 1990. "Multiple scales of statistical physics of neocortex: Application to electroencephalography." *Mathl. Comput. Modelling 13* no. 7: 83-95.

Ingber, L. and Nunez, P.L.. 1995. "Statistical mechanics of neocortical interactions: High resolution path-integral calculation of short-term memory." *Phys. Rev. E 51* no. 5: 5074-5083.

Ingber, L. and Rosen, B.. 1992. "Genetic algorithms and very fast simulated reannealing: A comparison." *Mathl. Comput. Modelling 16* no. 11: 87-100.

Ingber, L. and Sworder, D.D.. 1991. "Statistical mechanics of combat with human factors." *Mathl. Comput. Modelling 15* no. 11: 99-127.

Leonhard, R.. 1991. *The Art of Maneuver*. Presidio Press, Navato, CA.

Mayer, D.G., Belward, J.A., and Burrage, K.. 1996. "Use of advanced techniques to optimize a multi-dimensional dairy model." *Agricultural Systems 50*: 239-253.

Rosen, B.. 1992. "Function optimization based on advanced simulated annealing." *IEEE Workshop on Physics and Computation - PhysComp '92* : 289-293.

Sakata, S.. 1995. "High breakdown point estimation in econometrics". Ph.D. Thesis. University of California at San Diego, La Jolla, CA.

Sen, M.K. and Stoffa, P.L.. 1995. *Global Optimization Methods in Geophysical Inversion*. Elsevier, The NetherLands.

Su, A., Mager, S., Mayo, S.L., and Lester, H.A.. 1996. "A multi-substrate single-file model for ion-coupled transporters." *Biophys. J 70*: 762-777.

Toffler, A. and Toffler, H.. 1993. *War and Anti-War: Survival at the Dawn of the 21st Century*. Little, Brown, Boston.

Tzu, Sun. 1963. *The Art of War*. Oxford University Press, New York.

U.S. Army Training and Doctrine Command. 1996. *Land combat in the 21st century*. U.S. Army TRADOC, Washington, DC.

Wehner, M.F. and Wolfer, W.G.. 1983. "Numerical evaluation of path-integral solutions to Fokker-Planck equations. I.." *Phys. Rev. A 27*: 2663-2670.

Wofsey, M.. 1993. "Technology: Shortcut tests validity of complicated formulas." *The Wall Street Journal 222* no. 60: B1.

Wu, K. and Levine, M.D.. 1993. "3-D object representation using parametric geons". TR-CIM-93-13. Center for Intelligent Machines, McGill University, Montreal, Canada.